\input amstex
\input Amstex-document.sty

\pageno 233

\def\shorttitle{Recent Progress in Mathematical Analysis of Vortex Sheets}
\def\shortauthor{Sijue Wu}

\topmatter %
\title\nofrills{\boldHuge Recent Progress in Mathematical Analysis of Vortex Sheets}
\endtitle

\author \Large Sijue Wu* \endauthor

\thanks *Department of Mathematics, University of
Maryland, College Park, MD 20742, USA.
E-mail: sijue\@math.umd.edu
\endthanks

\abstract\nofrills \centerline{\boldnormal Abstract}

\vskip 4.5mm

{\ninepoint We consider the motion of the interface
 separating two domains of the same fluid that moves with different
velocity along the tangential direction of the interface. We assume
 that the fluids occupying the
two domains are of constant densities that are
equal,  are inviscid, incompressible
and irrotational, and that the surface tension is zero. We discuss
 results on the existence and uniqueness of solutions for given data,
 the regularity of solutions, singularity formation and the nature of
 solutions after the  singularity formation time.

\vskip 4.5mm


}
\endabstract
\endtopmatter

\document

\baselineskip 4.5mm \parindent 8mm

\specialhead \noindent \boldLARGE 1. Introduction \endspecialhead

Vortex dynamics is of fundamental importance for a wide variety of
concrete physical problems, such as lift of airfoils, mixing of
fluids, separation of boundary layers, and generation of sounds. In
mathematical analysis, one often neglects surface tension and
viscosity, when they are small in the real physical problem. This
necessitates justifying such simplifications.

In this paper, we consider the motion of the interface
separating two domains of the same fluid in $R^2$ that moves with
different velocity along the tangential direction of the interface.
We assume that the fluids occupying the
two domains separated by the interface are of constant densities that are
equal,  are inviscid, incompressible
and irrotational.
We also assume that the surface tension is zero, and there is no
external forces.  The
interface in the aforementioned fluid motion is a so called vortex sheet.
We want to study the following problem:

Given a vortex sheet initial data,
is there a unique solution to this problem?

In general, there are two
approaches to the aforementioned problem.
One is to solve the initial value problem of the incompressible
Euler equation in $R^2$:
$$\cases  v_t + v\cdot \nabla v +\nabla p=0
\text{div}\,v=0,\\v(x,y, 0)=v_0(x,y)\endcases\qquad (x,y)\in R^2, \ t\ge 0 \tag1$$ where the initial
incompressible velocity  $v_0\in L^2_{loc}(R^2)$, in which the vorticity $\omega_0=\text{curl}\,v_0$ is a finite
Radon measure. Here  $v$ is the fluid velocity, $p$ is the pressure, and  the density of the fluid is assumed to
be one. Notice that a vortex sheet gives a measure valued vorticity supported on the interface. This approach was
posed by DiPerna and Majda  in 1987 \cite{9}. In 1991, J.M\. Delort \cite{8} proved the existence of weak
solutions global in time of the 2-D incompressible Euler equation (1) for measure-valued initial vorticity in
$H^{-1}_{loc}(R^2)$ that has a distinguished sign.  However the problem of uniqueness of the weak solution is
still unresolved. In 1963, Yudovich \cite{33} obtained the existence and uniqueness of weak solutions of the 2-D
incompressible Euler equation (1) for initially bounded vorticity. The best results on uniqueness upto date are
given by Yudovich \cite{34} and Vishik \cite{31} for weak solutions with vorticity in a class slightly larger than
$L^\infty$. This does not include vortex sheets, which admit measure-valued vorticity. Examples of weak solutions
with the velocity field $v\in L^2( R^2\times (-T, T))$ that is compactly supported in space-time was constructed
by V. Scheffer \cite{29} and later by A. Shnirelman \cite{30}. This gives non-uniqueness of weak solutions in
$L^2(R^2\times (-T, T))$. However non-uniqueness in the physically relevant class of conserved energy $v\in
L^\infty([0,\infty),L^2_{loc}(R^2))$ remains open. Numerical evidences of non-uniqueness of weak solutions for
vortex sheet data can be found in \cite{25}, \cite{18}.

Furthermore, weak solutions give little information of the specific
nature of the
vortex sheet evolution. For instance, does the vorticity  remain supported
on a curve for a later time given that the initial vorticity is supported
on a curve in $R^2$? Assume further that the free interface between the
two fluid domains remains a curve in $R^2$ at
a later time,  equation (1) can be reduced to  an evolutionary
differential-integral equation along the
interface. This is the Birkhoff-Rott
equation,  written explicitly
 by Birkhoff in \cite{2} and implied in the work of Rott \cite{26}.
The second approach uses the Birkhoff-Rott equation as a model for the
evolution of the vortex sheet.

\specialhead \noindent \boldLARGE 2. The Birkhoff-Rott equation \endspecialhead

For convenience, we use complex variable $z=x+i\,y$ to
denote a point in $R^2$. $\overline{z}=x-i\,y$ denotes the complex
conjugate and
 $f_x=\partial_x f$ is the partial derivative of the function
$f$. $H^s$ indicates Sobolev spaces.

In search of the equation for the evolution of the vortex sheet, we
suppose that at time $t\ge 0$ the vorticity is a measure supported on the curve
$\Gamma(t)$ given by the complex position $\xi=\xi(s,t)$ in the arclength $s$,
in which $\xi(0,t)$ is the particle path of a
reference particle; and on this curve the
vorticity density is $\gamma=\gamma(s, t)$. That is the vorticity
at time $t$ is $\omega(x,y,t)$ satisfying
$$\iint \phi(x,y)\omega(x,y,t)\,dx\,dy=\int
\phi(\xi(s,t))\gamma(s,t)\,ds,\qquad \text{for any }\phi\in C_0^\infty(R^2).$$
From the Biot-Savart law, the velocity field $v$ induced
by the vorticity is given by
$$\overline{v}(z, t)=\frac1{2\pi
i}\int\frac{\gamma(s',t)}{z-\xi(s',t)}\,ds',\qquad\text{for }
z\not\in\Gamma(t).$$
Notice that the velocity
is discontinuous just on $\Gamma(t)$. We define the velocity on the
sheet as the average of the velocities at the two sides of the sheet,
that is given by the
principle value integral:
$$\overline{v}(\xi(s,t), t)=\frac1{2\pi i}\,p.v.\int\frac{\gamma(s',t)}{\xi(s,t)-\xi(s',t)}\,ds'.\tag2$$
As suggested by the properties of the Euler equation, we assume that
the vortex sheet is convected  by the
average velocity (2), and the vorticity
is conserved along the particle path. We arrive at the evolution equation of
the vortex sheet:
$$\aligned
\xi_t(s,t)+a(s,t)\xi_s(s,t)& =v(\xi(s,t),t)\\
\gamma_t(s,t)+\partial_s(a(s,t)\gamma(s,t))&=0\endaligned\tag3$$
where $a(s,t)$ is a real valued function satisfying $a(0,t)=0$.

A rigorous justification of the equivalence between  equation (3)
and equation (1) for smooth graphs $\xi(s,t)$ and smooth vortex
strength $\gamma(s,t)$ can be found in \cite{19}. It is not hard to
extend it to all smooth curves.

Assume $\alpha(s,t)=\int^s_0\gamma(s',t)\,ds'$
defines an increasing function of $s$, we make a change of variables:
$z(\alpha,
t)=\xi(s(\alpha, t),t)$,
in which $s(\alpha, t)$ is
the inverse of $\alpha(s,t)$: $\alpha(s(\alpha, t),t)=\alpha$.
We get from equation (3) the
 Birkhoff-Rott equation
$$
\partial_t\overline{z}(\alpha, t)=\frac 1{2\pi i}\,p.v.\int \frac
1{z(\alpha, t)-z(\beta, t)}\,d\beta.\tag4
$$
Notice that $z=z(\alpha, t)$ is a parameterization of the vortex sheet
in the circulation variable $\alpha$, $1/|z_\alpha|=\gamma$ is the vortex
strength. A steady solution of (4) is the flat sheet $z=\alpha$.

Equation (4) has been under active investigations over the last four
decades. A well-known property of (4) is that
perturbations of the flat sheet grow due to the Kelvin-Helmholtz
instability, following from a linearization of equation (4) about the
flat sheet.   For given analytic data, Sulem, Sulem,
Bardos and Frisch \cite{28} established the
short time existence and
uniqueness of solutions in analytic class for 2-D and 3-D vortex sheet
evolution. Duchon and Robert \cite{10} obtained the global
existence of solutions of equation (4) for a special class of initial data that
is close to the flat sheet. However,
numerous results show that a
vortex sheet can develop a curvature singularity in finite time from
analytic data. D.W. Moore \cite{21} was the first to provide
analytical evidence
that predicts the occurrence and time of singularity formation,  which
was verified numerically by Meiron, Baker and Orszag \cite{20} and by
Krasny \cite{13}. Caflisch and Orellana \cite{3} proved existence almost
up to the time of expected singularity formation for analytic data
that is close to the flat sheet. Duchon and Robert \cite{10} and
Caflisch and Orellana \cite{4} constructed specific examples
of solutions of equation (4) where a curvature singularity develops in
finite time from analytic data. The example of Caflisch and Orellana
\cite{4} has the form $z(\alpha, t)=\alpha+S(\alpha, t)+r(\alpha,t)$,
where
$$S(\alpha,t)=\epsilon(1-i)\{(1-e^{-t/2-i\alpha})^{1+\mu}-(1-e^{-t/2+i\alpha})^{1+\mu}\}$$
is a solution of the linearized equation in which $\epsilon$ is
small, $\mu>0$; $r(\alpha,t)$ is the correction term
that is negligible relative to $S(\alpha,t)$ in the
sense that $S(\alpha, t)+r(\alpha, t)$ exhibits the same kind of behavior as
$S(\alpha, t)$ \cite{4}. Notice that $S(\alpha,t)$ is
an analytic function for $t>0$, but $S(\alpha,0)$  has an infinite
second derivative at $\alpha=0$ for $\mu\in (0,1)$.
In fact, the $(1+\nu)$th derivative of $S(\alpha,0)$ for $\nu>\mu$
becomes infinite at $\alpha=0$.
Now inverting time gives an example $\hat{z}(\alpha, t)$ that is
analytic at $t_0<0$, but has an infinite second derivative at
$\alpha=0$, $t=0$. At the singularity formation time $t=0$, the vortex
strength $1/|\hat{z}_\alpha|$ of this example satisfies
$$0<c\le 1/|\hat{z}_\alpha| \le C<\infty\tag5$$
for some constants $c$ and $C$; and $\hat{z}(\alpha, t)\in
C^{1+\rho}(R\times[t_0, 0])$
for $0<\rho<\mu$.

These examples also show that the initial value problem of the Birkhoff-Rott
equation (4) is ill-posed in $C^{1+\nu}(R)$, $\nu>0$, and in Sobolev spaces
$H^s(R)$, $s>3/2$ in the Hadamard
sense
\cite{4}, \cite{10}. Ill-posedness was also
proved by Ebin \cite{11} using a  different approach. However the
existence of solutions in spaces less regular than $C^{1+\nu}(R)$ or
$H^s(R)$, and  the
nature of the vortex sheet at and beyond the singularity time
remained unknown analytically in general.

This suggests that we look for solutions of the
Birkhoff-Rott equation in the largest possible spaces where the
equation makes sense.
For the purpose of this paper, we consider
functions $z(\alpha, t)$
so that for each fixed time $t$, both sides of the equation (4) are
functions locally in $L^2$, and on
which the $L^2$-analysis is available. This leads us to consider
chord-arc curves,
thanks to the work of G. David \cite{7}.

Another reason that chord-arc curves are to be considered is due to
the numerical calculation of Krasny \cite{14} \cite{15}.
Krasny studied the evolution of the vortex sheet beyond singularity using
the vortex blob method.
He found that the approximating solutions
have the form of a spiral beyond
singularity. Convergence of the approximating sequence
to a weak solution of the Euler equation (1) was proved by J-G. Liu and
Z-P. Xin \cite{17}, under the assumption that the initial
vorticity has a distinguished sign. A special example of chord-arc
curves is a logarithmic spiral.


\specialhead \noindent \boldLARGE 3. Chord-arc curves and some recent
results  \endspecialhead

Let $\Gamma$  be a rectifiable Jordan curve in $R^2$ given by
$\xi=\xi(s)$ in the arclength $s$. We say $\Gamma$ is a chord-arc
curve, if there is a constant $M\ge 1$, such
that
$$|s_1-s_2|\le M|\xi(s_1)-\xi(s_2)|,\qquad\text{for all }s_1,\ s_2.$$
The infimum of all such constants M is called the chord-arc constant.

For a chord-arc curve $\xi=\xi(s)$, $s$ the arclength, it is proved in
\cite{6} that $\xi'(s)$ exists almost everywhere,
and there is a choice of the argument function $b\in BMO$, with
$\xi'(s)=e^{i b(s)}$. In particular, if the chord-arc constant is
close to 1, there is a choice of $b\in BMO$, such that
$\|b\|_{BMO}$ is close to $0$.
 Moreover the subset of all those functions $b$
is an open subset of $BMO$. And if $b\in BMO$, and
$\|b\|_{BMO}<1$,  $\xi(s)=\xi_0+\int_0^s e^{i\, b(s')}\,ds'$
defines a chord-arc curve.

Examples of chord-arc curves include Lipschitz curves and logarithmic
spirals $r=\pm e^\theta$, $\theta\in R$, where $(r,\theta)$ is the polar
coordinates.

A Theorem of G. David \cite{7} states that
\proclaim{Theorem (G. David \cite{7})} For all chord-arc curves $\Gamma:\,
\xi=\xi(s)$, $s$ the arclength,  the corresponding Cauchy integral
operator $C_\Gamma$, where
 $$C_\Gamma f(s)=p.v. \int \frac{f(s')}{\xi(s)-\xi(s')}\,d\,\xi(s'),$$
is bounded from $L^2(ds)$ to $L^2(ds)$.
\endproclaim

In fact, the result of G. David \cite{7} is stronger than stated
above. He proved that the Cauchy integral operator $C_\Gamma$ is
bounded from $L^2(ds)$ to $L^2(ds)$ if and only if $\Gamma$ is a regular
curve. A rectifiable curve $\Gamma$ is said to be regular if there is
a constant $M$ such that for every $r>0$ and every disc $D$ with radius
$r$, the length of $\Gamma\cap D$ does not exceed $Mr$. A chord-arc
curve is regular but not vice versa.

Now we go back to the Birkhoff-Rott equation. Notice that  the
Biot-Savart integral representing an incompressible velocity field $v$
in terms of the
vorticity $\omega$  may be divergent if
$\omega$ does not vanish fast enough at infinity, even if the
velocity field $v$ is well defined, we  extend the definition of
the Birkhoff-Rott equation (4) by considering the differences of the
velocities between any two points:
$$ \aligned
\overline{z_t}(\alpha, t)-&\overline{z_t}(\alpha', t)=\frac 1{2\pi
i}\,p.v.\int_{|\beta|\le N} \{\frac
1{z(\alpha, t)-z(\beta,
t)}-\frac
1{z(\alpha', t)-z(\beta,
t)}\}\,d\beta\\&+\frac 1{2\pi
i}\,p.v.\int_{|\beta|> N} \frac
{z(\alpha', t)-z(\alpha, t)}
{(z(\alpha, t)-z(\beta,
t))(z(\alpha', t)-z(\beta,
t))}\,d\beta,
\endaligned\tag6
$$
for all $(\alpha, t)$, $(\alpha',t)$,  and some
$N>|\alpha|+|\alpha'|+1$.
This  admits a larger
class of solutions. In particular, the integral on the right hand side of
equation (6) is  convergent for those similarity solutions
considered in \cite{12}, \cite{23}-\cite{25}, which otherwise give
divergent Cauchy integrals
in (4) due to the divergent contributions from the vorticities at infinity.
It follows from the Theorem of G. David that the integral on the right
hand side of the equation (6) is convergent for
$a.e.$ $(\alpha, t)$, $(\alpha',t)$ and is in $L^\infty([0, T],
L^2_{loc}(d\alpha)\times
L^2_{loc}(d\alpha'))$  for the solutions considered in Theorem 1 in the
following.

Roughly speaking, if a function $z=z(\alpha, t)$ satisfies equation
(4), it will also satisfy equation (6). On the other hand, if
$z=z(\alpha, t)$ satisfies (6), and the Cauchy
integral
$$\frac 1{2\pi i}\,p.v.\int \frac
1{z(\alpha, t)-z(\beta,
t)}\,d\beta$$
is convergent, then there is a function $c=c(t)$, such that $z(\alpha,
t)+c(t)$ satisfies equation (4).

For a local integrable function $f=f(\alpha)$ defined on $(a,b)$,
we say $f$ is of bounded local mean oscillation on $(a,b)$ if there exists
$\delta_0>0$ such that
$$\|f\|_{BMO(a,b), \delta_0}=\sup_{\text{all }I\subset (a,b),
|I|\le\delta_0}\frac 1{|I|}\int_I|f(\alpha)-f_I|\,d\alpha<\infty,$$
here $f_I=\frac 1{|I|}\int_I f(\alpha)\,d\alpha$, $I$ is an
interval. We say $f$ is
analytic on $(a, b)$ if $f\in C^\infty(a,b)$, and for any compact
subset $K$ of $(a, b)$,
there is a constant $\rho>0$, such that
$$\sum_{m=o}^{\infty}\frac{\rho^m}{m!}\int_K|\partial_\alpha^m
f(\alpha)|^2\,d\alpha<\infty.$$
Notice that $\ln z$ is multi-valued for complex number $z$. In the
following, $\ln z_\alpha$ refers to one choice of the multi-values.
We have the following results concerning the solutions of the
Birkhoff-Rott equation (6) (or (4)).

\proclaim{Theorem 1 \cite{32}} Assume that
$z\in H^1([0, T], L^2_{loc}(R))\cap L^2([0, T], H^1_{loc}(R))$ is a
solution of the
Birkhoff-Rott equation (6) for $0\le t\le T$, satisfying that

1. There are constants $m>0$, $M>0$, independent of $t$, such that
$$m|\alpha-\beta|\le |z(\alpha,t)-z(\beta,t)|\le M
|\alpha-\beta|\qquad\text{for all }\alpha,\ \beta,\ 0\le t\le T.\tag7$$

2. For the interval $(a,b)$, there exists $\delta_0>0$, independent of
   $t$, such that  for all $0\le t\le
   T$, $$\|\ln z_\alpha(\cdot,
   t)\|_{BMO(a,b), \delta_0}\le C(m, M),\tag8$$
here $C(m, M)$ is a
   universal constant depending on $m$ and $M$.
Assume further that  $ \ln z_\alpha\in L^2([0, T], L^2_{loc}(R))$.
Then  $ z_\alpha\in C((a,b)\times (0, T))$, and for
each $t_0\in (0, T)$, $z_\alpha(\cdot, t_0)$ is analytic
on $ (a,b)$.

\endproclaim

\noindent {\bf Remark } Notice that assumption 2. is satisfied if
$ \ln z_\alpha\in C([a, b]\times [0,T])$. Assumption 1.  is
equivalent to assuming that the vortex strength
$\gamma=1/|z_\alpha|$ is bounded away from $0$ and $\infty$, with
bounds independent of $t$, and $z=z(\cdot, t)$ defines a
chord-arc curve for each fixed $t\in [0,T]$, with the chord-arc
constant independent of $t$. Assumption 1. can be relaxed by
requiring that $1/|z_\alpha(\alpha, t)|$ is bounded away from $0$
and $\infty$ for $\alpha\in (a,b)$ only, and by assuming some
weaker conditions for $\alpha\not\in (a,b)$.

A version of Theorem 1 under assumptions that the solution of the
Birkhoff-Rott equation (4):
$z=z(\alpha, t)\in C^{1+\rho_0}$, for some $\rho_0>0$, the vortex strength
$c_0<\gamma(\alpha, t)<C_0$ for some constants $c_0>0$, $C_0>0$, and
the vortex sheets $\Gamma(t)$ are closed Jordan curves is also obtained by
Lebeau \cite{16} using an independent approach.

As a consequence of Theorem 1, the example constructed by Caflisch and
Orellana \cite{4} will fail to satisfy properties 1 and  2  as
stated in Theorem 1 near the singularity after the singularity formation time.

Let $\text{Re} z$ and $\text{Im} z$
be the real and imaginary parts of the complex number $z$ respectively.
Regarding the existence of solutions, we have
the following




\proclaim{Theorem 2 \cite{32}} For any real valued function $w_0\in H^{\frac
32}(R)$, there exists
$T=T(\,\|w_0\|_{H^{\frac32}}\,)>0$,
such that  the
Birkhoff-Rott equation (6) has a solution $z=z(\alpha, t)$ for $0\le
t\le T$, satisfying $\ln z_\alpha\in L^\infty([0,T],
H^{\frac32}(R))\cap Lip([0,T],
H^{\frac12}(R))$ and $\text{Im}\{ (1+i)\ln z_\alpha(\alpha, 0)\}=w_0(\alpha)$,
with the
properties that there exist constants $m>0$, $M>0$ independent of $t$,
such that
$$m|\alpha-\beta|\le |z(\alpha,t)-z(\beta,t)|\le M
|\alpha-\beta|,\qquad\text{for all }\alpha,\ \beta,\ 0\le t\le T;$$
 and there exists $\delta_0>0$, independent of $t$,
such that
$$\|\ln z_\alpha(\cdot,
   t)\|_{BMO(R), \delta_0}\le C(m, M),\qquad\text{for }0\le t\le T,$$
here $C(m, M)$ is the universal constant as in Theorem 1.
\endproclaim
Theorem 2
states that if only half of the data $z_\alpha(\alpha, 0)$ is given,
there is a solution of the Birkhoff-Rott equation (6) for a finite
time period. Theorem 2 is a generalization of the existence result of
Duchon and Robert \cite{10} to general data.

The following result implies that Theorem 2 is optimal, in the sense
that in general, there is no solution of the Birkhoff-Rott equation
satisfying properties 1. and 2. as stated in Theorem 1 beyond the initial time
$t=0$ for arbitrarily given data.

\proclaim{Theorem 3 \cite{32}} Assume that
$z\in H^1([0, T], L^2_{loc}(R))\cap L^2([0, T],
H^1_{loc}(R))$ is a solution of  the
Birkhoff-Rott equation (6) for $0\le t\le T$, $T>0$, satisfying the
property 1. and property 2.  on
some interval $(a, b)$ as stated
in Theorem 1.
Assume further that
$ \ln z_\alpha\in L^2([0, T], L^2_{loc}(R))$, and
$w_0=\text{Im}\{ (1+i)\ln
z_\alpha(\cdot, 0)\}$ is analytic on $(a,b)$. Then  $z_\alpha\in
C((a,b)\times [0,T))$ and
$\text{Re} \{ (1+i)\ln
z_\alpha(\cdot, 0)\}$ is also analytic on $(a,b)$.
\endproclaim


\specialhead \noindent \boldLARGE 4. Open questions \endspecialhead

The reason that Theorem 1-3 holds is that under the assumptions 1. and
2. in Theorem 1, the Birkhoff-Rott equation (6) is of ``elliptic''
type on $(a,b)\times[0,T]$. It would be interesting to see whether the
assumption 2 that requires small local mean oscillation on $\ln
z_\alpha(\cdot, t)$ can be
removed, or to see whether one can construct a
non-smooth solution of equation (6) that violates (8). A good place to
start is to construct similarity solutions of the Birkhoff-Rott
equation. Similarity solutions are studied numerically in Pullen
\cite{23} \cite{25}, Pullen and Phillips \cite{24}, and analytically in Kambe
\cite{12}.  However, the
similarity solutions in \cite{12}, \cite{23}-\cite{25}
violate (7). It would
also be interesting to
relax the assumption 1. by considering vortex strength that is not
necessarily bounded away from $0$ and infinity on the interval $(a,b)$.

Theorem 1-3 implies that a solution of the Birkhoff-Rott equation is
either analytic, or doesn't satisfy properties 1. and 2. in Theorem 1.
And there are in general no
solutions of reasonable regularity (properties 1. and 2.) beyond the
initial time for an
arbitrarily given data. Are there solutions of the Birkhoff-Rott
equation in spaces of even less regularity?  Notice that we can define
the solutions of the Birkhoff-Rott equation in the distribution sense by
$$\partial_t(\int \overline{z}(\alpha,t)\eta(\alpha)\,d\alpha)=\frac1{4\pi
i} \iint\frac{\eta(\alpha)-\eta(\beta)}{z(\alpha, t)-z(\beta,
t)}\,d\alpha\,d\beta,\qquad\text{for all }\eta\in C_0^\infty.$$
This admits some even larger classes of solutions $z=z(\alpha, t)$.
Is there a solution in the distribution sense for the Birkhoff-Rott
equation for any given data? How is it related to the weak solution of
the Euler equation?

Theorem 1-3 might as well suggest that the vortex sheet in general fails
to be a curve beyond the initial time for general data. Therefore it becomes
interesting to study the vortex layers or considering the effects of
viscosity. Numerical analysis of the vortex layers can be found in
Baker and Shelley \cite{1} etc.

\widestnumber \key {[10]}

\specialhead \noindent \boldLARGE References \endspecialhead

\ref \key 1 \by G. R. Baker \& M. J. Shelley \paper \rm On the connection between thin vortex layers and vortex
sheets \jour \it J. Fluid. Mech. \vol \rm 215 \yr 1990 \pages 161--194
\endref

\ref \key 2 \by G. Birkhoff \paper \rm Helmholtz and Taylor
instability \jour \it Proc. Symp. Appl. Math. \vol \rm XII, AMS
\yr 1962 \pages 55--76
\endref

\ref \key 3 \by R.E. Caflisch \& O.F. Orellana \paper \rm
Long-time existence for a slightly perturbed vortex sheet \jour
{\it Comm. Pure Appl. Math} \vol \rm 39 \yr 1986 \pages 807--838
\endref

\ref \key 4 \by R.E. Caflisch \& O.F. Orellana \paper \rm
Singular solutions and ill-posedness for the evolution of vortex
sheets \jour {\it SIAM J. Math. Anal.} \vol \rm 20 \yr 1989
\pages 293--307
\endref


\ref
\key 5
\by R. Coifman \& Y. Meyer
\paper \rm Lavrentiev's curves and conformal mappings
\jour \it Institut Mittag-Leffler
\vol \rm Report no.5
\yr 1983
\endref

\ref
\key 6
\by G. David
\book \it Th\'ese de troisi\'eme cycle, Meth\'ematique
\publ  Universit\'e de Paris XI
\endref

\ref \key 7 \by G. David \paper \rm Op\'erateurs int\'egraux
singuliers sur certaines courbes du plan complexe \jour \it  Ann.
Sci. \'Ecole. Norm. Sup.(4) \vol \rm 17, no.1 \yr 1984 \pages
157--189
\endref

\ref \key 8 \by J.M. Delort \paper \rm Existence de nappes de
tourbillon en dimension deux \jour {\it J. AMS} \vol \rm 4 \yr
1991 \pages 553--586
\endref

\ref \key 9 \by R. Diperna \& A. Majda \paper \rm Concentrations
in regularizations for 2-D incompressible flow \jour {\it Comm.
Pure Appl. Math.} \vol \rm 40 \yr 1987 \pages 301--345
\endref


\ref \key 10 \by J. Duchon \& R. Robert \paper \rm Global vortex
sheet solutions of Euler equations in the plane \jour {\it J.
Diff. eqns.} \vol \rm 73 \yr 1988 \pages 215--224
\endref

\ref \key 11 \by D.G. Ebin \paper \rm Ill-posedness of the
Rayleigh-Taylor and Helmholtz problems for incompressible fluids
\jour \it   Comm. PDE. \vol \rm 13 \yr 1988 \pages 1265--1295
\endref


\ref \key 12 \by T. Kambe \paper \rm Spiral vortex solutions of
Birkhoff-Rott equation \jour \it   Physica D \vol \rm 37 \yr 1989
\pages 463--473
\endref

\ref \key 13 \by R. Krasny \paper \rm On singularity formation in
a vortex sheet by the point-vortex approximation \jour \it  J.
Fluid Mech. \vol \rm 167 \yr 1986 \pages 65--93
\endref

\ref \key 14 \by R. Krasny \paper \rm Desingularization of
periodic vortex sheet roll-up \jour \it  J. Comput. Phys. \vol
\rm 65 \yr 1986 \pages 292--313
\endref

\ref \key 15 \by R. Krasny \paper \rm Computing vortex sheet
motion \jour \it  Proc. ICM'90, Kyoto, Japan \vol \rm Vol. II \yr
1991 \pages 1573--1583
\endref

\ref \key 16 \by G. Lebeau \paper \rm R\'egularit\'e du
probl\'eme de Kelvin-Helmholtz pour l'\'equation d'Euler 2D \jour
\it  S\'eminaire: \'Equations aux D\'eriv\'ees Partielles,
2000-2001, \vol \rm Exp\. No\. II \yr 2001 \pages 12
\endref

\ref \key 17 \by J-G. Liu \& Z-P. Xin \paper \rm Convergence of
vortex methods for weak solutions to the 2D Euler equations with
vortex sheet data \jour \it   Comm. Pure Appl. Math. \vol \rm
XLVIII \yr 1995 \pages 611--628
\endref

\ref
\key 18
\by M.C. Lopes, J. Lowengrub, H.J. Nussenzveig Lopes \& Y-X. Zheng
\paper \rm Numerical evidence for nonuniqueness evolution for the 2D
incompressible Euler equations: a vortex sheet example
\jour \it preprint
\endref

\ref
\key 19
\by C. Machioro \& M. Pulvirenti
\book \it Mathematical theory of incompressible nonviscous fluids
\publ Springer
\yr 1994
\endref


\ref
\key 20
\by D.I. Meiron, G.R. Baker \& S.A. Orszag
\paper \rm Analytic structure of vortex sheet dynamics. Part I.
\jour \it  J. Fluid Mech.
\vol \rm 114
\yr 1982
\pages 283
\endref

\ref \key 21 \by D.W. Moore \paper \rm The spontaneous appearance
of a singularity in the shape of an evolving vortex sheet \jour
\it  Proc. Roy. Soc. London Ser. A \vol \rm 365 \yr 1979 \pages
105--119
\endref

\ref \key 22 \by D.W. Moore \paper \rm Numerical and analytical
aspects of Helmholtz instability \jour \it   Theoretical and
Applied Mechanics, Proc. XVI ICTAM eds. Niordson and Olhoff \publ
North-Holland \yr 1984 \pages 629--633
\endref

\ref \key 23 \by D.I. Pullin \paper \rm The large scale structure
of unsteady self-similar roll-up vortex sheets \jour \it J. Fluid
Mech. \vol \rm 88, part 3 \yr 1978 \pages 401--430
\endref

\ref \key 24 \by D.I. Pullin \& W.R.C. Phillips \paper \rm On a
generalization of Kaden's problem \jour \it J. Fluid Mech. \vol
\rm 104 \yr 1981 \pages 45--53
\endref

\ref \key 25 \by D.I. Pullin \paper \rm On similarity flows
containing two branched vortex sheets \jour \it  Mathematical
aspects of vortex dynamics, ed. R. Caflisch \publ SIAM \yr 1989
\pages 97--106
\endref

\ref \key 26 \by N. Rott \paper \rm Diffraction of a weak shock
with vortex generation \jour \it  J. Fluid Mech. \vol \rm 1 \yr
1956 \pages 111--128
\endref

\ref
\key 27
\by P.G. Saffman
\book \it Vortex dynamics
\publ Cambridge
\yr 1992
\endref

\ref \key 28 \by P. Sulem, C. Sulem, C. Bados\& U. Frisch \paper
\rm Finite time analyticity for the two and three dimensional
Kelvin-Helmholtz instability \jour \it   Comm. Math. Phys. \vol
\rm 80 \yr 1981 \pages 485--516
\endref

\ref \key 29 \by V. Scheffer \paper \rm An inviscid flow with
compact support in space-time \jour \it  J. Geom. Anal. \vol \rm 3
\yr 1993 \pages 343--401
\endref

\ref \key 30 \by A. Shnirelman \paper \rm On the non-uniqueness
of weak solutions of the Euler equations \jour \it  Comm. Pure
Appl. Math \vol \rm L \yr 1997 \pages 1261--1286
\endref

\ref \key 31 \by M. Vishik \paper \rm Incompressible flows of an
ideal fluid with vorticity in borderline spaces of Besov type
\jour \it  Ann. Sci. \'Ecole Norm. Sup.(4) \vol \rm 32 no.6 \yr
1999 \pages 769--812
\endref

\ref
\key 32
\by S. Wu
\paper {\rm Recent progress in mathematical analysis of 2-D vortex sheet}
\jour preprint
\endref

\ref \key 33 \by V. Yudovich \paper \rm Non-stationary flow of an
ideal incompressible liquid \jour \it  USSR Comp. Math. and Math.
Phys. (English transl.) \vol \rm 3 \yr 1963 \pages 1407--1457
\endref

\ref \key 34 \by V. Yudovich \paper \rm Uniqueness theorem for
the basic nonstationary problem in the dynamics of an ideal,
incompressible fluid \jour \it  Math. Res. Lett. \vol \rm 2 \yr
1995 \pages 27--38
\endref

\enddocument